\documentclass[11pt,draft]{amsart}

\usepackage{cite,amsmath,amssymb,amsthm,latexsym,a4}

\usepackage[cp1251]{inputenc}

\newtheorem{theorem}{Theorem}[section]

\newtheorem{proposition}[theorem]{Proposition}

\newtheorem{lemma}[theorem]{Lemma}

\newtheorem{corollary}[theorem]{Corollary}

\DeclareMathOperator{\con}{con}

\DeclareMathOperator{\Gr}{Gr}

\numberwithin{equation}{section}

\makeatletter

\renewcommand*\subjclass[2][2010]{\def\@subjclass{#2}\@ifundefined{subjclassname@#1}{\ClassWarning{\@classname}{Unknown edition (#1) of Mathematics Subject Classification; using '2010'.}}{\@xp\let\@xp\subjclassname\csname subjclassname@#1\endcsname}}

\makeatother

\begin{document}
\title[Cancellable elements of the lattice of epigroup varieties]{Cancellable elements\\
of the lattice of epigroup varieties}

\author{Dmitry~V.~Skokov}

\address{Ural Federal University, Institute of Natural Sciences and Mathematics, Lenina 51, 620083 Ekaterinburg, Russia}

\email{dmitry.skokov@gmail.com}

\date{}

\thanks{The work was partially supported by the Ministry of Education and Science of the Russian Federation (project 1.6018.2017/8.9) and by the Russian Foundation for Basic Research (grant No.\,17-01-00551).}

\begin{abstract}
We completely determine all commutative epi\-group varieties that are cancellable elements of the lattice \textbf{EPI} of all epigroup varieties. In particular, we verify that a commutative epigroup variety is a cancellable element of the lattice \textbf{EPI} if and only if it is a modular element of this lattice.
\end{abstract}

\keywords{Epigroup, variety, lattice of varieties, EPI, modular element, cancellable element}

\subjclass{Primary 20M07, secondary 08B15}

\maketitle

\section{Introduction and summary} \label{introduction}
A semigroup $S$ is called an \emph{epigroup} if, for any element $x$ of $S$, there is a natural $n$ such that $x^n$ is a \emph{group element} (this means that $x^n$ lies in some subgroup of~$S$). Extensive information about epigroups may be found in the fundamental work~\cite{Shevrin-94} by L.\,N.\,Shevrin and the survey~\cite{Shevrin-05} by the same author. The class of epigroups is very wide. In particular, it includes all periodic semigroups (because some power of each element in such a semigroup lies in some its finite cyclic subgroup) and all completely regular semigroups (in which all elements are group elements). 

It is natural to consider epigroups as unary semigroups, i.e., semigroups equipped with an additional unary operation. The unary operation on an epigroup is defined in the following way. For an element $x$ of a given epigroup, let $e_x$ be the unit element of the maximal subgroup $G$ that contains some power of $x$. It is known (see~\cite{Shevrin-94}, for instance) that $xe_x = e_xx$ and this element lies in $G$. We denote by $\overline{x}$ the element inverse to $xe_x$ in $G$. This element is called the \emph{pseudoinverse} of $x$. The mapping $x\longmapsto\overline x$ defines a unary operation on an epigroup. Throughout this paper, we consider epigroups as algebras with the operations of multiplication and pseudoinversion. This naturally leads to the concept of varieties of epigroups as algebras with these two operations. An idea to examine epigroups in the framework of the theory of varieties was promoted by L.\,N.\,Shevrin in~\cite{Shevrin-94} (see also~\cite{Shevrin-05}). 

It is well known and may be easily checked that in every periodic epigroup the operation of pseudoinversion may be expressed in terms of multiplication (see~\cite{Shevrin-94}, for instance). This means that periodic varieties of epigroups may be identified with periodic varieties of semigroups.

There are many articles devoted to the examination of the special elements of the lattice \textbf{SEM} of all semigroup varieties and its sublattices.  A valuable information about  it can be found in~\cite{Vernikov-15} (see also \cite[\S 14]{Shevrin-Vernikov-Volkov-09}). Recently the analogues of a number of this results for the lattice \textbf{EPI} of all epigroup varieties was published in~\cite{Skokov-15, Skokov-16, Shaprynskii-Skokov-Vernikov-16}. The present article continues these investigations.   
	
In the lattice theory, special elements of many different types are considered. We recall definitions of three types of elements that appear below. An element $x$ of a lattice $\langle L;\vee,\wedge\rangle$ is called \emph{neutral} if
$$(\forall y,z\in L)\quad (x\vee y)\wedge(y\vee z)\wedge(z\vee x)= (x\wedge y)\vee(y\wedge z)\vee(z\wedge x).$$
It is well known that an element $x$ is neutral if and only if, for all $y,z\in L$, the sublattice of $L$ generated by $x$, $y$ and $z$ is distributive (see~\cite[Theorem~254]{Gratzer-11}). Neutral elements play an important role in the general lattice theory (see~\cite[Section~III.2]{Gratzer-11}, for instance). An element $x \in L$ is called
\begin{align*}
&\text{\emph{modular} if}&&(\forall y,z\in L)\quad y\le z\longrightarrow(x\vee y)\wedge z=(x\wedge z)\vee y,\\
&\text{\emph{cancellable} if}&&(\forall y,z\in L)\quad x\vee y=x\vee z\ \&\ x\wedge y=x\wedge z\longrightarrow y=z.
\end{align*}
It is easy to see that any cancellable element is a modular one. A valuable information about modular and cancellable elements in abstract lattices can be found in~\cite{Seselja-Tepavcevic-01}, for instance.

Modular elements of \textbf{SEM} and \textbf{EPI} were examined in ~\cite{Jezek-McKenzie-93, Shaprynskii-12, Vernikov-07} and~\cite{Shaprynskii-Skokov-Vernikov-16}, respectively. In particular, commutative semigroup varieties that are modular elements in \textbf{SEM} are completely determined in~\cite[Theorem 3.1]{Vernikov-07}. And the analogue of this statement in \textbf{EPI} was proved in~\cite[Theorem 1.3]{Shaprynskii-Skokov-Vernikov-16}. Further, it is verified in~\cite{Gusev-Skokov-Vernikov-17+} that the properties of being modular and cancellable elements of
\textbf{SEM} are equivalent in the class of commutative semigroup varieties. Here we are going to check this fact for a commutative epigroup varieties in \textbf{EPI}.

To formulate the main result of the article, we need some notation. We denote by $F$ the free unary semigroup
over a countably infinite alphabet. As usual, elements of $F$ are called \emph{words}. Words unlike variables are written in bold. Two parts of an identity we connect by the symbol~$\approx$, while the symbol~$=$ denotes the equality relation on $F$. Note that a semigroup $S$ satisfies the identity system $\mathbf wx\approx x\mathbf w\approx\mathbf w$ where the variable $x$ does not occur in the word \textbf w if and only if $S$ contains a zero element~0 and all values of \textbf w in $S$ are equal to~0. We adopt the usual convention of writing $\mathbf w\approx 0$ as a short form of such a system and referring to the expression $\mathbf w\approx 0$ as to a single identity. We denote by \textbf T the trivial semigroup variety and by \textbf{SL} the variety of all semilattices.

The main result of the article is the following

\begin{theorem}
\label{main}
For a commutative epigroup variety $\mathbf V$, the following are equivalent:
\begin{itemize}
\item[\textup{a)}] $\mathbf V$ is a cancellable element of the lattice $\mathbf{EPI}$;
\item[\textup{b)}] $\mathbf V$ is a modular element of the lattice $\mathbf{EPI}$;
\item[\textup{c)}] $\mathbf{V=M\vee N}$ where $\mathbf M$ is one of the varieties $\mathbf T$ or $\mathbf{SL}$, while $\mathbf N$ is a variety satisfying the identities $x^2y\approx 0$ and $xy\approx yx$.
\end{itemize}
\end{theorem}

Theorem~\ref{main} and~\cite[Theorem~1.3]{Shaprynskii-Skokov-Vernikov-16} imply the following

\begin{corollary}
\label{cancellable and modular in SEM and EPI}
For a periodic commutative semigroup variety $\mathbf V$, the following are equivalent:
\begin{itemize}
\item[\textup{a)}]$\mathbf V$ is a cancellable element of the lattice $\mathbf{EPI}$;
\item[\textup{b)}]$\mathbf V$ is a modular element of the lattice $\mathbf{EPI}$;
\item[\textup{c)}]$\mathbf V$ is a cancellable element of the lattice  $\mathbf{SEM}$;
\item[\textup{d)}]$\mathbf V$  is a modular element of the lattice $\mathbf{SEM}$.\qed
\end{itemize}
\end{corollary}

The article consists of three sections. Section~\ref{prel} contains an auxiliary facts, while Section~\ref{proof} is devoted to verification of Theorem~\ref{main}.

\section{Preliminaries}
\label{prel}
The following two lemmas are well known and may be easily checked. (see~\cite{Shevrin-94,Shevrin-05}, for instance).
\begin{lemma}
\label{x^omega} If $S$ is an arbitrary epigroup and $x\in S$ then $x\,\overline{x}\,=x^\omega$.\qed
\end{lemma}

This Lemma allows us to use for brevity the expression of the kind $a^\omega$ where $a\in F$ as a short form of the word $a\,\overline{a}$.  The following is well known (see \cite{Shevrin-94,Shevrin-05}, for instance).

\begin{lemma}
\label{in every epigroup}
The identities
\begin{equation}
\label{eq in every epigroup} x^{\omega}x \approx xx^{\omega} \approx \overline{\overline{x}} \approx x^{\omega}\,\overline{\overline{x}} \approx \overline{\overline{x}}\,x^{\omega}
\end{equation}
hold in an arbitrary epigroup.\qed
\end{lemma}

The following assertion may be easily checked.

\begin{lemma}
\label{x*=0}
The identity $\overline x\approx 0$ holds in an epigroup $S$ if and only if $S$ is a nil-semigroup.\qed
\end{lemma}

For a given word \textbf w, we denote by $\con(\mathbf w)$ the \emph{content} of \textbf w, i.e., the set of all variables occurring in \textbf w. A word written in the language of multiplication and pseudoinversion is called a \emph{semigroup word} if it does not include the operation of pseudoinversion. If \textbf w is a semigroup word then $\ell(\mathbf w)$ stands for the length of \textbf w; otherwise, we put $\ell(\mathbf w) = \infty$. If nilpotency index of any nil-semigroup in \textbf X is not exceeded some natural number $n$ and $n$ is the least number with such a property then $n$ is called a \emph{degree} of the variety \textbf X and is denoted by $\deg(\mathbf X)$; otherwise we put $\deg(\mathbf X)=\infty$. The following proposition is verified in~\cite[Corollary 1.3]{Gusev-Vernikov-16} and plays an essential role in the further.
\begin{proposition}
	\label{deg in EPI}
	For an epigroup variety $\mathbf V$, the following are equivalent: 
		\begin{itemize}
		\item[\textup{(i)}] $\deg(\mathbf V)\le n$;
		\item[\textup{(ii)}] $\mathbf V$ satisfies an identity of the form  $x_1 x_2 \cdots x_n \approx \mathbf v$ where $\con(\mathbf v) = \{x_1,x_2, \dots, x_n\}$ and $\ell(\mathbf v) > n$;
		\item[\textup{(iii)}] $\mathbf V$ satisfies an identity of the form 
		\begin{equation}
		\label{degEPI}
		x_1 x_2 \cdots x_n \approx x_1 x_2 \cdots x_{i-1}\cdot \overline{\overline{x_i \cdots x_j}}\cdot x_{j+1}\cdots x_n
		\end{equation} for some $1 \le i \le j \le n$.\qed
	\end{itemize}
\end{proposition}

Proposition~\ref{deg in EPI} immediately implies the following
\begin{corollary}
\label{degree of meet}
The equality $\deg (\mathbf X \wedge \mathbf Y) = \min \bigl\{\deg (\mathbf X), \deg (\mathbf Y) \bigr\}$ holds in an arbitrary epigroup variety.\qed
\end{corollary}

Repeating literally the proof of Corollary 2.13 of~\cite{Vernikov-08-umod1}, but using Proposition~\ref{deg in EPI} of this paper instead of Proposition 2.1 of~\cite{Vernikov-08-umod1}, we see that the following statement is true.

\begin{corollary}
\label{degree of join with nil}
If $\mathbf V$~is an arbitrary epigroup variety and $\mathbf N$~is a nil-variety then $\deg (\mathbf V \vee \mathbf N) = \max\bigl\{ \deg(\mathbf V), \deg (\mathbf N)\bigr\}$.\qed
\end{corollary}

The first statement of the following lemma is generally known (see~\cite{Shevrin-94,Shevrin-05}). The second claim is verified  in~\cite[Theorem~1.1]{Shaprynskii-Skokov-Vernikov-16}.

\begin{lemma}
\label{SL is neutral atom}
The variety $\mathbf{SL}$ is
\begin{itemize}
\item[\textup{(i)}] an atom of the lattice $\mathbf{EPI}$;
\item[\textup{(ii)}]  a neutral element of $\mathbf{EPI}$.\qed
\end{itemize}
\end{lemma}

The following lemma is well known and may be easily checked.

\begin{lemma}
\label{word problem}
An identity $\mathbf{u}\approx\mathbf{v}$ holds in the variety $\mathbf {SL}$ if and only if $\con(\mathbf{u})=\con(\mathbf{v})$.\qed
\end{lemma}

We need the following three lemmas.

\begin{lemma}[{\mdseries\cite[Lemma~2.1]{Gusev-Skokov-Vernikov-17+}}]
\label{join with neutral atom} Let $L$ be a lattice with~$0$ and $a$ an atom and neutral element of $L$. An element $x \in L$ is cancellable if and only if the element $x \vee a$ is cancellable.\qed
\end{lemma}

\begin{lemma}[{\mdseries\cite[Lemma~2.2]{Gusev-Skokov-Vernikov-17+}}]
\label{over neutral atom}
Let $L$ be a lattice with~$0$, $a$ an atom and neutral element of $L$ and $x\in L$. If, for any $y,z\in L$, the equalities $x\vee(y\vee a)=x\vee(z\vee a)$ and $x\wedge(y\vee a)=x\wedge(z\vee a)$ imply that $y\vee a=z\vee a$ then $x$ is a cancellable element.\qed
\end{lemma}

\begin{lemma}[{\mdseries\cite[Lemma~2.3]{Gusev-Skokov-Vernikov-17+}}]
\label{modular non-cancellable}
Let $x$ be a modular but not cancellable element of a lattice $L$ and let $y$ and $z$ be different elements of $L$ such that $x\vee y=x\vee z$ and $x\wedge y=x\wedge z$. Then there is an element $x'\in L$ such that $x'\le x$, $x'\vee y=x'\vee z$, $x'\wedge y=x'\wedge z$ and $y\vee z=x'\vee y$.\qed
\end{lemma}

\section{Proof of the main result}
\label{proof}
In this section we prove Theorem~\ref{main}. The implication a)~$\longrightarrow$~b) is evident, while the equivalence of the claims b) and c) is checked in~\cite[Theorem~1.3]{Shaprynskii-Skokov-Vernikov-16}. It remains to prove the implication c)~$\longrightarrow$~a). Lemmas~\ref{SL is neutral atom} and~\ref{join with neutral atom} allow us to assume that $\mathbf V=\mathbf N$. Suppose that \textbf N is non-cancellable element of $\mathbf{EPI}$. Hence there are epigroup varieties \textbf Y and \textbf Z with $\mathbf{N\vee Y}=\mathbf{N\vee Z}$, $\mathbf{N\wedge Y}=\mathbf{N\wedge Z}$ and $\mathbf{Y\ne Z}$.

Repeating literally the proof of Lemma 3.1 of~\cite{Gusev-Skokov-Vernikov-17+}, but using the corollaries~\ref{degree of meet} and~\ref{degree of join with nil} of this paper instead of Lemmas 2.5 and 2.6 of~\cite{Gusev-Skokov-Vernikov-17+}, we see that the following statement is true.

\begin{lemma}
	\label{degree Y = degree Z} $\deg (\mathbf Y) = \deg (\mathbf Z).$\qed
\end{lemma}

Since the claims~b) and~c) of Theorem~\ref{main} are equivalent, \textbf N is a modular element of $\mathbf{EPI}$. In view of Lemma~\ref{modular non-cancellable}, there is a variety $\mathbf N'$ such that
$$
\mathbf{N'\subseteq N},\,\mathbf{N'\vee Y}=\mathbf{N'\vee Z}=\mathbf{Y\vee Z}\text{ and }\mathbf{N'\wedge Y}=\mathbf{N'\wedge Z}.
$$

Since $\mathbf{Y\ne Z}$, we can assume without loss of generality that there is an identity $\mathbf{u\approx v}$ that holds in \textbf Y but is false in \textbf Z. If this identity is satisfied by the variety $\mathbf N'$ then it holds in $\mathbf{N'\vee Y}=\mathbf{N'\vee Z}$, and therefore, in \textbf Z. Thus, $\mathbf{u\approx v}$ is wrong in $\mathbf N'$. By Lemma~\ref{x*=0}, one of the words \textbf u and \textbf v is not equal to 0 in $ \mathbf N'$. So we can assume without loss of generality that $\mathbf u$ is a semigroup word. Recall that a semigroup word \textbf w is called \emph{linear} if any variable occurs in \textbf w at most once. Being a subvariety of \textbf N, the variety $\mathbf N'$ satisfies the identities $x^2y\approx 0$ and $xy\approx yx$. Therefore, any non-linear word except $x^2$ is equal to~0 in $\mathbf N'$. Thus, we may assume without loss of generality that either $\mathbf u=x^2$ or $\mathbf u=x_1x_2\cdots x_k$ for some $k$. It is evident that the varieties $\mathbf Y $ and $\mathbf Z $ are either both periodic or both non-periodic, because otherwise one of the varieties $\mathbf N ' \vee \mathbf Y$ or $ \mathbf N' \vee \mathbf Z$ is periodic and the other one is not. If $\mathbf Y$ and $\mathbf Z$ are periodic we may consider it as a semigroup varieties. In this case the proof of Theorem~\ref{main} may be complete by repeating literally arguments from the proof of Theorem 1~\cite{Gusev-Skokov-Vernikov-17+}. So we may assume that $\mathbf Y$ and $\mathbf Z$ are non-periodic. Lemmas~\ref{SL is neutral atom} and~\ref{over neutral atom} allow us to assume that $\mathbf Y, \mathbf Z \supseteq \mathbf{SL}$. Lemma~\ref{word problem} implies that $\con(\mathbf u)=\con(\mathbf v)$. Combining the observations given above, we have that $\mathbf{u\approx v}$ is either an identity of the form  $x^2 \approx \mathbf w$ where $\mathbf w$~is not a semigroup word and $\con(\mathbf w) = \{x\}$ or $x_1x_2 \cdots x_k \approx \mathbf v$ where $\con(\mathbf v) = \{x_1, x_2, \dots, x_k\}$ and $\ell(\mathbf v) \ge k$.  If $\ell(\mathbf v) = k$ then the identity $x_1x_2 \cdots x_k \approx \mathbf v$ holds in $\mathbf N'$, a cotradiction. So $\ell(\mathbf v) > k$. 

Further considerations may be divided into two cases depending on the form of the identity $\mathbf{u\approx v}$.

\smallskip
\emph{Case 1}: $\mathbf{u\approx v}$ is an identity of the form $x^2 \approx \mathbf w$, where $\mathbf w$ is not a semigroup word and $\con(\mathbf w) = \{x\}$. Lemma~\ref{x*=0} implies that  $\mathbf N'$ satisfies the identity $\mathbf w \approx 0$. If $\mathbf N'$ satisfies $x^2\approx 0$ then the identity $x^2\approx \mathbf w$ holds in $\mathbf{N'\vee Y}=\mathbf{N'\vee Z}$. We see that the identity $x^2\approx \mathbf w$ holds in \textbf Z, a contradiction. So, we may assume that the identity $x^2\approx 0$ does not hold in $\textbf N'$. Since \textbf Y satisfies the identity $x^2\approx \mathbf w$, the same identity holds in $\mathbf{N'\wedge Y}=\mathbf{N'\wedge Z}$. Therefore, there is a deduction of this identity from identities of the varieties $\mathbf N'$ and \textbf Z. In particular, one of these varieties satisfies a non-trivial identity of the form  $x^2\approx\mathbf w'$. Suppose at first that this identity holds in $\mathbf N'$.  If $\mathbf w'$ is not a semigroup word that $\mathbf N'$ satisfies the identity $x^2\approx 0$ by Lemma~\ref{x*=0}. And we have a contradiction in this case. We may assume that  $\mathbf w'$ is a semigroup word. We substitute 0 into any symbol from $\con(\mathbf w')\setminus\{x\}$ whenever $\con(\mathbf w')\ne\{x\}$. We obtain that $x^2\approx 0$ holds in $\mathbf N'$ and we have a contradiction. If $\con(\mathbf w')=\{x\}$ then \textbf Y is a periodic variety, a contradiction. Therefore, the identity $x^2\approx\mathbf w'$ does not hold in $\mathbf N'$. This means that this identity holds in \textbf Z. Since $\mathbf{Z\supseteq SL}$, Lemma~\ref{word problem} implies that $\con (\mathbf w') = \{x\}$. Further, $\mathbf w'$ is not a semigroup word, while \textbf Z is a non-periodic variety. 

By~\cite[Lemma 2.4]{Gusev-Vernikov-16}, if \textbf w is not a semigroup word and $\con(\mathbf w)=\{x\}$ that an arbitrary epigroup variety satisfies the identity $\mathbf{w} \approx x^p{\,\overline{x}\,}^q$ for some integer $p \ge 0$ and some natural $q$. So, there are $p\ge 0$ and $q\in\mathbb N$ such that the identity $x^2 \approx x^p\,\overline{x}\,^q$ holds in $\mathbf Y$ but does not hold in $\mathbf Z$. We denote by $\Gr S$ the set of all group elements of an epigroup $S$. Let $S\in\mathbf Y$ and $x\in S$.

Suppose that $p\le q$. By~\ref{x^omega} we obtain that
$$
x^2=x^p\,\overline{x}\,^q=x^\omega\,\overline x\,^{q-p}\in \Gr S.
$$
Finally, let $p>q$. Lemma~\ref{x^omega} and the identities~\eqref{eq in every epigroup} imply that
$$
x^2=x^p\,\overline{x}\,^q=x^{p-q}x^\omega=x^{p-q}(x^\omega)^{p-q}=\,\overline{\overline x}\,^{^{p-q}}\in \Gr S.
$$
So, $x^2\in\Gr S$. Analogously, the identity $x^2\approx\mathbf w'$ holds in \textbf Z. This means that $x^2\in\Gr S$ always for any epigroup $S\in\mathbf Z$ and any $x\in S$. Clearly, the varieties $\mathbf Y$ and $\mathbf Z$ satisfy the identity $x^2 \approx x^2 x^{\omega}$. Lemmas ~\ref{x^omega} and~\ref{x*=0} imply that $\mathbf N'$ satisfies the identity
$$
x^2 x^{\omega}\approx x^3\,\overline x\approx 0\approx x^p\,\overline{x}\,^q.
$$
The identities $x^2 x^\omega\approx x^2\approx x^p\,\overline{x}\,^q$ also hold in $\mathbf Y$. So, the variety $\mathbf{N'\vee Y}=\mathbf{N'\vee Z}$ satisfies the identity $x^2 x^{\omega} \approx x^p\,\overline{x}\,^q$. Therefore, this identity holds in $\mathbf Z$. This means that \textbf Z satisfies the identities $x^2 \approx x^2 x^{\omega} \approx x^p\,\overline{x}\,^q$ and we have a contradiction.

\smallskip
\emph{Case 2}: $\mathbf{u\approx v}$ is an identity of the form $x_1x_2\cdots x_k\approx\mathbf v$ where $\con(\mathbf v)=\{x_1,x_2,\dots,x_k\}$ and $\ell(\mathbf v) > k$. Let $\deg (\mathbf Y) = n$. In view of Lemma~\ref{degree Y = degree Z}, $\deg (\mathbf Z) = \deg (\mathbf Y) = n$. Proposition~\ref{deg in EPI} implies that $n \le k$. Lemma~\ref{modular non-cancellable} allow us to assume that $\mathbf Y \vee \mathbf Z = \mathbf N' \vee \mathbf Y = \mathbf N' \vee \mathbf Z$. It is clear that $\deg (\mathbf Y \vee \mathbf Z) \ge n$. Suppose at first that $\deg (\mathbf Y \vee \mathbf Z) = n$. Then
$$
\deg(\mathbf N') \le \deg (\mathbf N' \vee \mathbf Y) = \deg (\mathbf Y \vee \mathbf Z) = n.
$$
Being a nil-variety, $\mathbf N'$ satisfies the identity $x_1x_2\cdots x_n\approx 0$ in this case. Since $\ell(\mathbf v)>k\ge n$, the identity $x_1x_2\cdots x_k\approx\mathbf v$ holds in $\mathbf N'$ as well. This contradicts the choice of the identity $\mathbf{u\approx v}$.

Let now $\deg (\mathbf Y \vee \mathbf Z) > n$. Since $\deg (\mathbf Y) = n$, Lemma~\ref{deg in EPI} implies that \textbf Y satisfies an identity of the form~\eqref{degEPI} for some $1 \le i \le j \le n$. The same lemma implies that this identity is false in $\mathbf{Y\vee Z}$.  Therefore,~\eqref{degEPI} is wrong in $\mathbf Z$. Analogously, there are  $1 \le i' \le j' \le n$ such that the identity
\begin{equation}
\label{degEPI'} x_1 x_2 \cdots x_n \approx x_1 x_2 \cdots x_{i'-1}\cdot \overline{\overline{x_{i'} \cdots x_{j'}}}\cdot x_{j'+1}\cdots x_n
\end{equation}
holds in \textbf Z but does not hold in \textbf Y. We will assume without any loss that $i \le i'$. Further considerations may be divided into three subcases, depending on the relations between $j$, $i '$, and $j'$.

\smallskip
\emph{Subcase 2.1}: Suppose at first that $i\le j < i' \le j'$. Then we substitute $({x_{i'}\cdots x_{j'}})^{\omega}x_{j'+1}$ into $x_{j'+1}$ in~\eqref{degEPI} whenever $j'<n$ or multiply~\eqref{degEPI} by $({x_{i'} \cdots x_{j'}})^{\omega}$ on the right whenever $j'=n$. As a result, we obtain the identity
$$
\begin{array}{rl}
 & x_1 \cdots x_i \cdots x_j \cdots (x_{i'}  \cdots x_{j'})(x_{i'} \cdots x_{j'})^{\omega}x_{j'+1} \cdots x_n \\[3pt]
\approx & x_1 \cdots \overline{\overline{x_i \cdots x_j}} \cdots x_{i'} \cdots x_{j'}(x_{i'} \cdots x_{j'})^{\omega}x_{j'+1} \cdots x_n.
\end{array}
$$
Clearly, this identity holds in $\mathbf N'$. So, it also holds in $\mathbf Z$ because $\mathbf N' \vee \mathbf Y = \mathbf N' \vee \mathbf Z$. This identity and the identity~\eqref{eq in every epigroup} imply that $\mathbf Z$ satisfies the identity
\begin{equation}
\label{3 no intersection}
x_1 \cdots x_i \cdots x_j \cdots \overline{\overline{x_{i'} \cdots x_{j'}}} \cdots x_n \approx x_1 \cdots \overline{\overline{x_i \cdots x_j}} \cdots \overline{\overline{x_{i'} \cdots x_{j'}}} \cdots x_n.
\end{equation}

Substitute $x_{i-1}(x_i \cdots x_j)^{\omega}$ into $x_{i-1}$ in~\eqref{degEPI'} whenever $i>1$ or multiply~\eqref{degEPI'} by $({x_{i} \cdots x_{j}})^{\omega}$ on the left whenever $i=1$. We obtain
$$
\begin{array}{rl}
&x_1 \cdots x_{i-1}(x_i  \cdots x_j)^{\omega}x_i \cdots x_{j} \cdots x_{i'} \cdots x_{j'} \cdots x_n  \\[3pt]
\approx & x_1 \cdots x_{i-1}(x_i \cdots x_j)^{\omega}x_i \cdots x_{j} \cdots  \overline{\overline{x_{i'} \cdots x_{j'}}} \cdots x_n.
\end{array}
$$
This identity and the identity~\eqref{eq in every epigroup} imply that $\mathbf Z$ satisfies the identity
\begin{equation}
\label{4 no intersection}
x_1 \cdots \overline{\overline{x_i  \cdots x_j}} \cdots x_{i'} \cdots x_{j'} \cdots x_n \approx x_1 \cdots \overline{\overline{x_i \cdots x_j}} \cdots \overline{\overline{x_{i'} \cdots x_{j'}}} \cdots x_n.
\end{equation}
And this means that \textbf Z satisfies the identities
\begin{align*}
x_1 x_2 \cdots x_n\stackrel{\eqref{degEPI'}}\approx{}& x_1 x_2 \cdots x_{i'-1}\cdot \overline{\overline{x_{i'} \cdots x_{j'}}}\cdot x_{j'+1}\cdots x_n\\[-3pt]
\stackrel{\eqref{3 no intersection}}\approx{}&x_1 \cdots \overline{\overline{x_i \cdots x_j}} \cdots \overline{\overline{x_{i'} \cdots x_{j'}}} \cdots x_n\\[-3pt]
\stackrel{\eqref{4 no intersection}}\approx{}&x_1 \cdots \overline{\overline{x_i  \cdots x_j}} \cdots x_{i'} \cdots x_{j'} \cdots x_n.
\end{align*}
Here we write $\mathbf{w\stackrel{\varepsilon}\approx w'}$ in the case when the identity $\mathbf{w\approx w'}$ follows from the identity $\varepsilon$. So, we prove that the identity~\eqref{degEPI} holds in \textbf Z, a contradiction. 

\smallskip
\emph{Subcase 2.2}: Suppose that $i < i' \le j < j'$. Substitute $({x_{i'}\cdots x_j \cdots x_{j'}})^{\omega}x_{j'+1}$ into $x_{j'+1}$ in~\eqref{degEPI} whenever $j'<n$ or multiply~\eqref{degEPI} by $({x_{i'}\cdots x_j \cdots x_{j'}})^{\omega}$ on the right whenever $j'=n$. We obtain the identity
\begin{equation}
\label{3}
\begin{array}{rl}
 & x_1 \cdots x_i \cdots (x_{i'} \cdots x_j \cdots x_{j'})(x_{i'}\cdots x_j \cdots x_{j'})^{\omega}x_{j'+1} \cdots x_n \\[3pt]
\approx & x_1 \cdots \overline{\overline{x_i \cdots x_{i'} \cdots x_j}}\cdots x_{j'}(x_{i'}\cdots x_j \cdots x_{j'})^{\omega}x_{j'+1} \cdots x_n.
\end{array}
\end{equation}
It is evident that this identity holds in $\mathbf N'$. The equality $\mathbf N' \vee \mathbf Y = \mathbf N' \vee \mathbf Z$ implies that it holds in \textbf Z too. Substitute $x_{i-1}(x_i \cdots x_{i'} \cdots x_j)^{\omega}$ into $x_{i-1}$ in~\eqref{degEPI'} whenever $i>1$ or multiply~\eqref{degEPI'} by $({x_{i}\cdots x_{i'} \cdots x_{j}})^{\omega}$ on the left whenever $i=1$. We obtain the identity
\begin{equation}
\label{4}
\begin{array}{rl}
&x_1 \cdots x_{i-1}(x_i \cdots x_{i'} \cdots x_j)^{\omega}x_i\cdots x_{i'}\cdots x_{j} \cdots x_{j'} \cdots x_n  \\[3pt]
\approx & x_1 \cdots x_{i-1}(x_i \cdots x_{i'} \cdots x_j)^{\omega}x_i \cdots  \overline{\overline{x_{i'}\cdots x_j \cdots x_{j'}}} \cdots x_n.
\end{array}
\end{equation}
Clearly, this identity holds in $\mathbf Z$. As a result we obtain that $\mathbf Z$ satisfies the identities
\begin{align*}
x_1x_2\cdots x_n \stackrel{\eqref{degEPI'}}\approx{}&x_1 \cdots x_{i} \cdots \overline{\overline{x_i' \cdots x_j \cdots x_j'}}\cdot x_{j'+1}\cdots x_n\\[-4pt]
\stackrel{\eqref{eq in every epigroup}}\approx{}&x_1 \cdots x_i \cdots (x_{i'} \cdots x_j \cdots x_{j'})(x_{i'}\cdots x_j \cdots x_{j'})^{\omega}x_{j'+1} \cdots x_n\\[-3pt]
\stackrel{\eqref{3}}\approx{}&x_1 \cdots \overline{\overline{x_i \cdots x_{i'} \cdots x_j}}\cdots x_{j'}(x_{i'}\cdots x_j \cdots x_{j'})^{\omega}x_{j'+1} \cdots x_n\\[-3pt]
\stackrel{\eqref{eq in every epigroup}}\approx{}&x_1 \cdots (x_i \cdots x_j)^{\omega}(x_i \cdots x_j)\cdots x_{j'}(x_{i'}\cdots x_j \cdots x_{j'})^{\omega}\cdot x_{j'+1} \cdots x_n\\[-3pt]
\stackrel{\eqref{eq in every epigroup}}\approx{}&x_1 \cdots x_{i-1}(x_i \cdots x_j)^{\omega}x_i \cdots  x_{i'-1}\cdot \overline{\overline{x_{i'}\cdots x_j \cdots x_{j'}}}\cdot x_{j'+1} \cdots x_n\\[-3pt]
\stackrel{\eqref{4}}\approx{}&x_1 \cdots x_{i-1}(x_i \cdots x_j)^{\omega}x_i\cdots x_{i'}\cdots x_{j} \cdots x_{j'} \cdots x_n\\[-3pt]
\stackrel{\eqref{eq in every epigroup}}\approx{}&x_1 x_2 \cdots x_{i-1}\cdot \overline{\overline{x_i \cdots x_j}}\cdot x_{j+1}\cdots x_n.
\end{align*}
So, $\mathbf Z$ satisfies~\eqref{degEPI}. We have a contradiction. 

\smallskip

\emph{Subcase 2.3}: Finally suppose that $i \le i' \le j' \le j$. Clearly, $\mathbf Z$ satisfies the identity~\eqref{degEPI} whenever $i = i'$ and $j = j'$. We have a contradiction in this case. Therefore, we will assume without any loss that $i < i'$. Substitute $x_{i'-1}(x_{i'}\cdots x_{j'})^{\omega}$ into $x_{i'-1}$ in~\eqref{degEPI}. By the identity~\eqref{eq in every epigroup}, we obtain that
\begin{equation}\label{23x}
x_1 \cdots x_{i'-1}\cdot \overline{\overline{x_{i'} \cdots x_{j'}}}\cdot x_{j'+1}\cdots x_n \approx x_1 \cdots\overline{\overline{ x_{i}\cdots \overline{\overline{x_{i'} \cdots x_{j'}}}\cdots x_{j}}}\cdots x_n.
\end{equation}
Substitute $x_{i-1}(x_i \cdots \overline{\overline{x_{i'} \cdots x_{j'}}}\cdots x_{j})^{\omega}$ into $x_{i-1}$ in~\eqref{degEPI'} whenever $i>1$ or multiply~\eqref{degEPI'} by $(x_i \cdots \overline{\overline{x_{i'} \cdots x_{j'}}}\cdots x_{j})^{\omega}$ on the left whenever $i=1$. By Lemma~\ref{x^omega}, we obtain the identity
\begin{equation}\label{33x}
\begin{array}{rl}
&x_1  \cdots x_{i-1}(x_i \cdots \overline{\overline{x_{i'} \cdots x_{j'}}}\cdots x_{j})^{\omega} x_i \cdots x_{i'} \cdots x_{j'} \cdots x_j \cdots x_n  \\[2pt]
\approx &x_1 \cdots\overline{\overline{ x_{i}\cdots \overline{\overline{x_{i'} \cdots x_{j'}}}\cdots x_{j}}}\cdots x_n.
\end{array}
\end{equation}
So, $\mathbf Z$ satisfies the identities
\begin{align*}
x_1x_2\cdots x_n \stackrel{\eqref{degEPI'}}\approx{}&x_1 \cdots x_i \cdots x_{i'-1}\cdot \overline{\overline{x_{i'} \cdots x_{j'}}}\cdot x_{j'+1}\cdots x_n\\[-3pt]
\stackrel{\eqref{23x}}\approx{}&x_1 \cdots\overline{\overline{ x_{i}\cdots \overline{\overline{x_{i'} \cdots x_{j'}}}\cdots x_{j}}}\cdots x_n\\[-3pt]
\stackrel{\eqref{33x}}\approx{}&x_1 \cdots (x_{i}\cdots \overline{\overline{x_{i'} \cdots x_{j'}}}\cdots x_{j})^{\omega}(x_{i}\cdots x_{i'} \cdots x_{j'}\cdots x_{j})\cdots x_n.
\end{align*}
We see that the identities
\begin{equation}
\label{one more} x_1x_2\cdots x_n \approx x_1 \cdots\overline{\overline{ x_{i}\cdots \overline{\overline{x_{i'} \cdots x_{j'}}}\cdots x_{j}}}\cdots x_n
\end{equation}
and
\begin{equation}
\label{two more} x_1x_2\cdots x_n \approx x_1 \cdots (x_{i}\cdots \overline{\overline{x_{i'} \cdots x_{j'}}}\cdots x_{j})^{\omega}(x_{i}\cdots x_{i'} \cdots x_{j'}\cdots x_{j})\cdots x_n
\end{equation}
hold in $\mathbf Z$.
Substitute $(x_i \cdots x_{i'} \cdots x_{j'} \cdots x_j)^{\omega}x_{j+1}$ into $x_{j+1}$ in~\eqref{two more} whenever $j<n$ or multiply this identity by $(x_i \cdots x_{i'} \cdots x_{j'} \cdots x_j)^{\omega}$ on the right whenever $j=n$. We obtain the identity
\begin{equation}
\label{*}
\begin{array}{rl}
&x_1 \cdots \overline{\overline{x_i \cdots x_{i'} \cdots x_{j'} \cdots x_j}} \cdots x_n \\[2pt]
\approx & x_1 \cdots x_{i-1}(x_i \cdots \overline{\overline{x_{i'} \cdots x_{j'}}} \cdots x_j)^{\omega}\cdot \overline{\overline{x_i \cdots x_{i'} \cdots x_{j'} \cdots x_j}}\cdots x_n.
\end{array}
\end{equation}
Then substitute $x_{i-1}(x_i \cdots \overline{\overline{x_{i'} \cdots x_{j'}}} \cdots x_j)^{\omega}$ into $x_{i-1}$ in~\eqref{degEPI} whenever $i>1$ or multiply~\eqref{degEPI} by $(x_i \cdots \overline{\overline{x_{i'} \cdots x_{j'}}} \cdots x_j)^{\omega}$ whenever $i=1$. We obtain the identity
\begin{equation}
\label{22}
\begin{array}{rl}
&x_1 \cdots x_{i-1}(x_i \cdots \overline{\overline{x_{i'} \cdots x_{j'}}} \cdots x_j)^{\omega}x_i \cdots x_{i'} \cdots x_{j'} \cdots x_j \cdots x_n \\[2pt]
\approx & x_1 \cdots x_{i-1}(x_i \cdots \overline{\overline{x_{i'} \cdots x_{j'}}} \cdots x_j)^{\omega}\cdot \overline{\overline{x_i \cdots x_{i'} \cdots x_{j'} \cdots x_j}}\cdots x_n.
\end{array}
\end{equation}
By the same way, we substitute $x_{i-1}(x_i \cdots \overline{\overline{x_{i'} \cdots x_{j'}}} \cdots x_j)^{\omega}$ into $x_{i-1}$ in~\eqref{one more} whenever $i>1$ or multiply~\eqref{one more} by $(x_i \cdots \overline{\overline{x_{i'} \cdots x_{j'}}} \cdots x_j)^{\omega}$ on the left whenever $i=1$. We obtain the identity
\begin{equation}
\label{23}
\begin{array}{rl}
& x_1 \cdots x_{i-1}(x_i \cdots \overline{\overline{x_{i'} \cdots x_{j'}}} \cdots x_j)^{\omega}x_i \cdots x_{i'} \cdots x_{j'} \cdots x_j \cdots x_n \\[3pt]
\approx & x_1 \cdots x_{i-1}(x_i \cdots \overline{\overline{x_{i'} \cdots x_{j'}}} \cdots x_j)^{\omega}\cdot \overline{\overline{x_i \cdots \overline{\overline{x_{i'} \cdots x_{j'}}} \cdots x_j}}\cdots x_n.
\end{array}
\end{equation}
Finally, $\mathbf Z$ satisfies the identities
\begin{align*}
x_1x_2\cdots x_n \stackrel{\eqref{degEPI'}}\approx{}&x_1 x \cdots x_i \cdots \overline{\overline{x_{i'} \cdots x_{j'}}}\cdots x_j \cdots x_n\\[-3pt]
\stackrel{\eqref{23x}}\approx{}&x_1  \cdots\overline{\overline{ x_{i}\cdots \overline{\overline{x_{i'} \cdots x_{j'}}}\cdots x_{j}}}\cdots x_n\\[-3pt]
\stackrel{\eqref{eq in every epigroup}}\approx{}&x_1 \cdots x_{i-1}(x_i \cdots \overline{\overline{x_{i'} \cdots x_{j'}}} \cdots x_j)^{\omega}\cdot \overline{\overline{x_i \cdots \overline{\overline{x_{i'} \cdots x_{j'}}} \cdots x_j}}\cdots x_n\\[-3pt]
\stackrel{\eqref{23}}\approx{}&x_1 \cdots x_{i-1}(x_i \cdots \overline{\overline{x_{i'} \cdots x_{j'}}} \cdots x_j)^{\omega}x_i \cdots x_{i'} \cdots x_{j'} \cdots x_j \cdots x_n \\[-3pt]
\stackrel{\eqref{22}}\approx{}&x_1 \cdots x_{i-1}(x_i \cdots \overline{\overline{x_{i'} \cdots x_{j'}}} \cdots x_j)^{\omega}\cdot \overline{\overline{x_i \cdots x_{i'} \cdots x_{j'} \cdots x_j}}\cdots x_n\\[-3pt]
\stackrel{\eqref{*}}\approx{}&x_1 \cdots \overline{\overline{x_i \cdots x_{i'} \cdots x_{j'} \cdots x_j}} \cdots x_n  .
\end{align*}
We prove that the identity~\eqref{degEPI} holds in \textbf Z, a contradiction. This completes the proof of Theorem~\ref{main}.\qed

\medskip

\textbf {Acknowledgements}. The author would like to thank Professor B. M. Vernikov for helpful discussions and several useful remarks.


\label{lastpage}  

\end{document}